\documentstyle[11pt]{article}
\begin{document}
\newcommand{\2}{\vspace{0.2 cm}}
\newcommand{\dom}{\mbox{$\rightarrow$}}
\newcommand{\ndom}{\mbox{$\not\rightarrow$}}
\newcommand{\sdom}{\mbox{$\Rightarrow$}}
\newcommand{\nsdom}{\mbox{$\not\Rightarrow$}}
\newcommand{\qed}{\hfill$\Box$}
\newcommand{\pf}{{\bf Proof: }}
\newtheorem{theorem}{Theorem}[section]
\newtheorem{algorithm}[theorem]{Algorithm}
\newtheorem{proposition}[theorem]{Proposition}
\newtheorem{lemma}[theorem]{Lemma}
\newtheorem{problem}[theorem]{Problem}
\newtheorem{corollary}[theorem]{Corollary}
\newtheorem{conjecture}[theorem]{Conjecture}
\newcommand{\beq}{\begin{equation}}
\newcommand{\eeq}{\end{equation}}
\newcommand{\ra}{\rangle}
\newcommand{\la}{\langle}
\newcommand{\har}{\rightleftharpoons}
\newcommand{\<}[1]{\mbox{$\la #1 \ra$}}

\title{Hamilton Cycles in Digraphs of Unitary Matrices}

\date{}

\author{
G. Gutin\thanks{Department of Computer Science, Royal Holloway,
University of London, Egham, Surrey, TW20 0EX, UK,
gutin@cs.rhul.ac.uk}\and  A. Rafiey\thanks{Department of Computer
Science, Royal Holloway, University of London, Egham, Surrey, TW20
0EX, UK, arash@cs.rhul.ac.uk} \and S. Severini\thanks{Department of
Mathematics and Department of Computer Science, University of York,
York, YO10 5DD, UK, ss54@york.ac.uk} \and A. Yeo\thanks{Department
of Computer Science, Royal Holloway, University of London, Egham,
Surrey, TW20 0EX, UK, anders@cs.rhul.ac.uk}}

\maketitle

\begin{abstract}
A set $S\subseteq V$ is called an {\em $q^+$-set} ({\em $q^-$-set},
respectively) if $S$ has at least two vertices and, for every $u\in
S$, there exists $v\in S, v\neq u$ such that $N^+(u)\cap N^+(v)\neq
\emptyset$ ($N^-(u)\cap N^-(v)\neq \emptyset$, respectively). A
digraph $D$ is called {\em s-quadrangular} if, for every $q^+$-set
$S$, we have $|\cup \{N^+(u)\cap N^+(v): u\neq v, u,v\in S\}|\ge
|S|$ and, for every $q^-$-set $S$, we have $|\cup \{N^-(u)\cap
N^-(v): u,v\in S)\}\ge |S|$. We conjecture that every strong
s-quadrangular digraph has a Hamilton cycle and provide some support
for this conjecture.

\end{abstract}

\section{Introduction}

The hamiltonian cycle problem is one of the central problems in
graph theory and its applications
\cite{bang2000,gouldGC19,west1996}. Many sufficient conditions were
obtained for hamiltonicity of undirected graphs \cite{gouldGC19} and
only a few such conditions are proved for directed graphs (for
results and conjectures on sufficient conditions for hamiltonicity
of digraphs, see \cite{bang2000}). This indicates that the asymmetry
of the directed case makes the Hamilton cycle problem significantly
harder, in a sense.

For a digraph $D=(V,A)$ and $x\neq y\in V$, we say that $x$ {\em
dominates} $y$, denoted $x\dom y$, if $xy\in A.$ All vertices
dominated by $x$ are called the {\em out-neighbors} of $x$; we
denote the set of out-neighbors by $N^+(x)$. All vertices that
dominate $x$ are {\em in-neighbors} of $x$; the set of in-neighbors
is denoted by $N^-(x).$ A set $S\subseteq V$ is called an {\em
$q^+$-set} ({\em $q^-$-set}, respectively) if $S$ has at least two
vertices and, for every $u\in S$, there exists $v\in S, v\neq u$
such that $N^+(u)\cap N^+(v)\neq \emptyset$ ($N^-(u)\cap N^-(v)\neq
\emptyset$, respectively). A digraph $D$ is called {\em
s-quadrangular} if, for every $q^+$-set $S$, we have $|\cup
\{N^+(u)\cap N^+(v): u\neq v, u,v\in S\}|\ge |S|$ and, for every
$q^-$-set $S$, we have $|\cup \{N^-(u)\cap N^-(v): u,v\in S)\}|\ge
|S|$. A digraph $D$ is {\em strong} if there is a path from $x$ to
$y$ for every ordered pair $x,y$ of vertices of $D.$

We believe that the following claim holds:

\begin{conjecture} \label{conj} Every strong s-quadrangular digraph is hamiltonian.
\end{conjecture}

A complex $n\times n$ matrix $U$ is \emph{unitary} if $U\cdot
U^{\dagger }=U^{\dagger}\cdot U=I_{n}$, where $U^{\dagger}$ denotes
the conjugate transpose of $U$ and $I_{n}$ the $n\times n$ identity
matrix. The \emph{digraph of }an $n\times n$ matrix $M$ (over any
field) is a digraph on $n$ vertices with an arc $ij$ if and only if
the $(i,j)$-entry of the $M$ is nonzero.  It was shown in
\cite{severiniSIAMJMAA25} that the digraph of a unitary matrix is
s-quadrangular; s-quadrangular tournaments were studied in
\cite{lundgren}.

It follows that if Conjecture \ref{conj} is true, then the digraph
of an irreducible unitary matrix is hamiltonian. Unitary matrices
are important in quantum mechanics and, at present, are central in
the theory of quantum computation \cite{nielsen2000}. In particular,
we may associate a strong digraph to a quantum system whose unitary
evolution allows transitions only along the arcs of the digraph
(that is, respecting the topology of the graph, like in discrete
quantum walks \cite{aharonov2001,kempeCP44}). Then, if the
conjecture is true, the digraph would be necessarily hamiltonian.
Moreover, the conjecture is important in the attempt to understand
the combinatorics of unitary and unistochastic matrices, see, e.g.,
\cite{beasleyPIMA50,fiedler1988,zyczkowskiJP36}. If the conjecture
is true, then the digraph of an irreducible weighing matrix has a
Hamilton cycle (see \cite{geramita1979}, for a reference on weighing
matrices). Also, since the Kronecker product of unitary matrices
preserves unitarity, if $K$ and $H$ are digraphs of irreducible
unitary matrices, then their Kronecker product $K \otimes H$ (see
\cite{imrich2000}, for an interesting collection of notions and
results on graph products) has a Hamilton cycle provided $K \otimes
H$ is strong. The {\em complete biorientation} of an undirected
graph $G$ is a digraph obtained from $G$ by replacing every edge
$xy$ by the pair $xy,yx$ of arcs. A graph is {\em s-quadrangular} if
its complete biorientation is s-quadrangular. Certainly, the
following is a weakening of Conjecture \ref{conj}:

\begin{conjecture} \label{conj1} Every connected s-quadrangular graph is hamiltonian.
\end{conjecture}

In this paper, we provide some support to the conjectures. In
Section \ref{di}, we show that if a strong s-quadrangular digraph
$D$ has the maximum semi-degree at most 3, then $D$ is hamiltonian.
In our experience, to improve the result by replacing
$\Delta^0(D)\le 3$ with $\Delta^0(D)\le 4$ appears to be a very
difficult task. In Section \ref{undi}, we show the improved result
only for the case of undirected graphs. Even in this special case
the proof is fairly non-trivial. Before recalling some standard
definitions and proving our results, it is worth mentioning that the
line digraphs of eulerian digraphs are all s-quadrangular and
hamiltonian. We have verified Conjecture \ref{conj} for all digraphs
with at most five vertices and a number of digraphs with six
vertices.

The number of out-neighbors (in-neighbors) of $x$ is the {\em
out-degree} $d^+(x)$ of $x$ ({\em in-degree} $d^-(x)$ of $x$). The
maximum {\em semi-degree} $\Delta^0(D)=\max\{ d^+(x),d^-(x):\ x\in
V\}.$ A collection of disjoint cycles that include all vertices of
$D$ is called a {\em cycle factor} of $D$. We denote a cycle factor
as the union of cycles $C_1\cup \cdots \cup C_t$, where the cycles
$C_i$ are disjoint and every vertex of $D$ belongs to a cycle $C_j.$
If $t=1$, then clearly $C_1$ is a {\em Hamilton cycle} of $D$. A
digraph with a Hamilton cycle is called {\em hamiltonian}. Clearly,
the existence of a cycle factor is a necessary condition for a
digraph to be hamiltonian.

\section{Supporting Conjecture \ref{conj}}\label{di}

The existence of a cycle factor is a natural necessary condition for
a digraph to be hamiltonian \cite{gutinJCT78}. The following
necessary and sufficient conditions for the existence of a cycle
factor is well known, see, e.g., Proposition 3.11.6 in
\cite{bang2000}.

\begin{lemma}\label{lem}
A digraph $H$ has a cycle factor if and only if, for every
$X\subseteq V(H)$, $|\cup _{x\in X}N^+(x)|\ge |X|$ and $|\cup _{x\in
X}N^-(x)|\ge |X|.$
\end{lemma}

Using this lemma, it is not difficult to prove the following
theorem:

\begin{theorem}\label{cfdi}
Every strong s-quadrangular digraph $D=(V,A)$ has a cycle factor.
\end{theorem}

\pf Let $X\subseteq V$. If $X$ is a $q^+$-set, then $$|X|\le |\cup
(N^+(u)\cap N^+(v): u\neq v, u,v\in X)|\le |\cup _{x\in X}N^+(x)|.$$
If $X$ is not a $q^+$-set, then consider a maximal subset $S$ of
$X$, which is a $q^+$-set (possibly $S=\emptyset$). Since $D$ is
strong every vertex of $X$ dominates a vertex. Moreover, since every
vertex of $X-S$ dominates a vertex that is not dominated by another
vertex in $X$, we have $|\cup_{x\in X-S}N^+(x)|\ge |X-S|.$ Thus,
$$|\cup _{x\in X}N^+(x)|\ge |X-S|+|\cup _{x\in S}N^+(x)|\ge
|X-S|+|S|=|X|.$$ Similarly, we can show that $|\cup _{x\in
X}N^-(x)|\ge |X|$ for each $X\subseteq V.$ Thus, by Lemma \ref{lem},
$D$ has a cycle factor. \qed

\vspace{2mm}

Now we are ready to prove the main result of this section.

\begin{theorem} If the out-degree and in-degree of every vertex in a strong
s-quadrangular digraph $D$ are at most 3, then $D$ is hamiltonian.
\end{theorem}
\pf Suppose that $D=(V,A)$ is a non-hamiltonian strong s-quadrangular digraph and for every
vertex $u\in V$, $d^+(u),d^-(u)\le 3.$

Let $F=C_1\cup\cdots \cup C_t$ be a cycle factor of $D$ with minimum
number $t\ge 2$ of cycles. Assume there is no cycle factor $C'_1\cup
\cdots\cup C'_t$ such that $|V(C'_1)| < |V(C_1)|$. For a vertex $u$
on $C_i$ we denote by $u^+$ $(u^-)$ the successor (the predecessor)
of $u$ on $C_i$. Also, define $x^{++}=(x^+)^+$. Since every vertex
belongs to exactly one cycle of $F$ these notations define unique
vertices. Let $u,v$ be vertices of $C_i$ and $C_j$, $i\neq j,$
respectively and let $K(u,v)=\{uv^+,vu^+\}$. At least one of the
arcs in $K(u,v)$ is not in $D$ as otherwise we may replace the pair
$C_i,C_j$ of cycles in $F$ with just one cycle $uv^+v^{++}\ldots
vu^+u^{++}\ldots u$, a contradiction to minimality of $t.$

Since $D$ is strong, there is a vertex $x$ on $C_1$ that dominates a
vertex $y$ outside $C_1$. Without loss of generality, we may assume
that $y$ is on $C_2$. Clearly, $\{x,y^-\}$ is a $q^+$-set. Since
$K(x,y^-)\not\subset A$ and $x\dom y$, we have $y^- \ndom x^+.$
Since $D$ is s-quadrangular, this is impossible unless $d^+(x)>2.$
So, $d^+(x)=3$ and there is a vertex $z\not\in \{x^+,y\}$ dominated
by both $x$ and $y^-.$ Let $z$ be on $C_j.$

Suppose that $j\neq 1.$ Since $K(x,z^-) \not\subset A$ and $x\dom
z$, we have $z^-\ndom x^+$. Since $\{x,z^-\}$ is a $q^+$-set, we
have $z^-\dom y$. Suppose that $j\ge 3.$ Since $z^-\dom y$ and
$y^-\dom z$, we have $K(y^-,z^-)\subset A,$ which is impossible. So,
$j=2$. Observe that $\{x^+,z\}$ is a $q^-$-set ($x\dom x^+$ and
$x\dom z$). But $y^-\ndom x^+$ and $z^-\ndom x^+$. Hence,
$|N^-(x^+)\cap N^-(z)|=1$, which is impossible.

Thus, $j=1.$ Since $K(y^-,z^-)\not\subset A$ and $y^-\dom z$, we
have $z^-\ndom y.$ Since $\{x,z^-\}$ is a $q^+$-set, $z^-\dom x^+.$
By replacing $C_1$ and $C_2$ in $F$ with $x^+x^{++}\ldots z^-x^+$
and $xyy^+\ldots y^-zz^+\ldots x$ we get a cycle factor of $D$, in
which the first cycle is shorter than $C_1$. This is impossible by
the choice of $F.$\qed

\section{Supporting Conjecture \ref{conj1}}\label{undi}

Let $G=(V,E)$ be an undirected graph and let $f:\ V \dom
\mathcal{N}$ be a function, where $\mathcal{N}$ is the set of
positive integers. A spanning subgraph $H$ of $G$ is called an
$f$-{\em factor} if the degree of a vertex $x\in V(H)$ is equal to
$f(x).$ Let $e(X,Y)$ denote the number of edges with one endpoint in
$X$ and one endpoint in $Y$. We write $e(X)=e(X,X)$ to denote the
number of edges in the subgraph $G\<{X}$ of $G$ induced by $X$. The
number of neighbors of a vertex $x$ in $G$ is called the {\em
degree} of $x$ and it is denoted by $d_G(x).$

The following assertion is the well-known Tutte's $f$-factor Theorem
(see, e.g., Exercise 3.3.16 in \cite{west1996}):

\begin{theorem} A graph $G=(V,E)$ has an $f$-factor if and only if
$$q(S,T)+\sum_{t\in T}(f(t)-d_{G-S}(t))\le \sum_{s\in S}f(s)$$ for
all choices of disjoint subsets $S,T$ of $V$, where $q(S,T)$ denotes
the number of components $Q$ of $G-(S\cup T)$ such that
$e(V(Q),T)+\sum_{v\in V(Q)}f(v)$ is odd.
\end{theorem}

The following lemma is of interest for arbitrary undirected graphs.
A 2-factor of $G$ is an {\em $f$-factor} such that $f(x)=2$ for each
vertex $x$ in $G.$

\begin{lemma} \label{NoCycFac}
If $G=(V,E)$ is a graph of minimum degree at least 2 and with no
2-factor, then we can partition $V(G)$ into $S$, $T$, $O$ and $R$,
such that the following properties hold.
\begin{description}
 \item[(i):] $T$ is independent.
 \item[(ii):] $e(R,O \cup T)=0$.
 \item[(iii):] Every connected component in $G\<{O}$ has an odd
       number of edges into $T$.
 \item[(iv):] No $t \in T$ has two edges into the same connected component
       of $G\<{O}$.
 \item[(v):] For every vertex $o \in O$ we have $e(o,T) \leq 1$.
 \item[(vi):] There is no edge $ot \in E(G)$, where $t \in T$, $o \in O$, such that $e(t,S)=0$
       and $e(o,O)=0$.
 \item[(vii):] $|T|-|S|-\frac{e(T,O)-oc(S,T)}{2} >0$, where $oc(S,T)$ is the
number of connected components in $G-S-T$, which have an odd
number of edges into $T$. (Note
          that $oc(S,T)$ is also the number of connected components
          of $G\<{O}$, by (ii) and (iii).)
\end{description}
\end{lemma}

\pf By Tutte's $f$-factor Theorem, there exists disjoint subsets $S$
and $T$ of $V$, such that the following holds.
$$ oc(S,T)+2|T|-\sum_{v \in T} d_{G-S}(v) > 2|S| $$
Define
 $$ w(S,T) =|T|-|S|-e(T)-\frac{e(T,V-S-T)-oc(S,T)}{2}. $$
We now choose disjoint subsets $S$ and $T$ of $V$, such that the
following holds in the order it is stated.
\begin{itemize}
  \item maximize $w(S,T)$
  \item minimize $|T|$
  \item maximize $|S|$
  \item minimize $oc(S,T)$
\end{itemize}
Furthermore, let $O$ contain all vertices from $V-S-T$ belonging to
connected components of $G\<{V-S-T}$, each of which has an odd
number of edges into $T$. Let $R=V-S-T-O$. We will prove that
$S,T,O$ and $R$ satisfy (i)-(vii).

Clearly, $w(S,T)>0$. Let $t \in T$ be arbitrary and assume that
$t$ has edges into $i$ connected components in $G\<{O}$. Let
$S'=S$ and let $T'=T-\{t\}$. Furthermore let  $j=1$ if the
connected component in $G\<{V-S'-T'}$ containing $t$ has an odd
number of edges into $T'$, and $j=0$, otherwise. Observe that
$|T'|=|T|-1$, $e(T')=e(T)-e(t,T)$, $e(T',V-S'-T')=e(T,V-S-T)-
e(t,V-S-T)+e(t,T)$ and $oc(S',T')=oc(S,T)-i+j$. Since $|T'|<|T|$,
we must have $w'=w(S',T')-w(S,T)<0.$ Therefore, we have

$$
w'=-1+e(t,T)+\frac{e(t,V-S-T)-e(t,T)-i+j}{2} <0.$$ Thus,
\begin{equation}\label{in1} e(t,T)+e(t,V-S-T)-i+j \le 1\end{equation}

Observe that, by the definition of $i$, $e(t,V-S-T)\ge i.$ Now, by
(\ref{in1}), $i\le e(t,V-S-T)\le i-j-e(t,T)+1.$ Thus, $j+e(t,T)\le
1$ and, if $e(t,T)>0$, then $e(t,T)=1$, $j=0$ and $e(t,V-S-T)=i$.
However, $e(t,T)=1$, $e(t,V-S-T)=i$ and a simple parity argument
imply $j=1$, a contradiction.

So $e(t,T)=0$. In this case, $e(t,V-S-T)-i=0$ or 1. If
$e(t,V-S-T)-i=1$, then, by a simple parity argument, we get $j=1$,
which is a contradiction against (\ref{in1}). Therefore, we must
have $e(t,V-S-T)=i$.

It follows from $e(t,T)=0$, $e(t,V-S-T)=i$, $w(S,T)>0$ and the
definition of $O$ that (i), (ii), (iii), (iv) and (vii) hold. We
will now prove that (v) and (vi) also hold. Suppose that there is
a vertex $o \in O$ with $e(o,T)>1$. Let $S'=S \cup \{o\}$ and
$T'=T$, and observe that
\begin{equation}\label{in2}
w(S',T')-w(S,T)=-1+(e(o,T)-(oc(S,T)-oc(S',T')))/2<0\end{equation} If
$e(o,T)\ge 3$ then, by (\ref{in2}), $oc(S,T)> oc(S',T')+1.$ However,
by taking $o$ from $O$, we may decrease $oc(S,T)$ by at most 1, a
contradiction. So, $e(o,T)=2$ and, by (\ref{in2}), $oc(S,T)>
oc(S',T').$ However, since $e(o,T)$ is even, taking $o$ from $O$
will not decrease $o(S,T)$, a contradiction. Therefore (v) holds.

Suppose that there is an edge $ot \in E(G)$, where $t \in T$, $o
\in O$, such that $e(t,S)=0$ and $e(o,O)=0$.   Let $S'=S$ and
$T'=T\cup \{o\} -\{t\}$. By (iii) and (iv), the connected
component in $G\<{V-S'-T'}$, which contains $t$, has an odd number
of arcs into $T'$. This implies that
$oc(S',T')-oc(S,T)=-e(t,O)+1$. Also,
$e(T,V-S-T)=e(T',V-S'-T')+e(t,O)-1.$ By (v) we have $e(T')=0$. The
above equalities imply that $w(S',T')=w(S,T)$. Since the degree of
$t$ is at least 2 and $e(t,S)=0$, we conclude that $e(t,O)\ge 2$.
Thus, $oc(S',T') < oc(S,T)$, which is a contradiction against the
minimality of $oc(S,T)$. This completes that proof of (vi) and
that of the lemma. \qed

For a vertex $x$ in a graph $G$, $N(x)$ denotes the set of neighbors
of $x$; for a subset of $X$ of $V(G)$, $N(X)=\cup_{x\in X}N(x).$ A
2-factor contains no cycle of length 2. Thus, the following theorem
cannot be deduced from Theorem \ref{cfdi}. For a vertex $x$ in a
graph $G$, $N(x)$ is the set of neighbors of $x$.

\begin{theorem}\label{2-factor}
Every connected s-quadrangular graph with at least three vertices
contains a 2-factor.
\end{theorem}

\pf Suppose $G$ is a connected $s$-quadrangular graph with at least
three vertices and with no 2-factor. Suppose $G$ has a vertex $x$ of
degree 1 such that $xy$ is the only edge incident to $x$. Consider
$z\neq x$ adjacent with $y$. Observe that vertices $x,z$ only one
common neighbor. Thus, $G$ is not s-quadrangular. Thus, we may
assume that the minimum degree of a vertex in $G$ is at least two
and we can use Lemma \ref{NoCycFac}.

Let $S$, $T$, $O$ and $R$ be defined as in Lemma \ref{NoCycFac}.
First suppose that there exists a vertex $t \in T$ with $e(t,S)=0$.
Let $y \in N(t)$ be arbitrary, and observe that $y \in O$, by (i)
and (ii). Furthermore observe that $e(y,O) \geq 1$ by (vi), and let
$z \in N(y) \cap O$ be arbitrary. Since $y \in N(t) \cap N(z)$,
there must exist a vertex $u \in N(t) \cap N(z) -\{y\}$, by the
definition of an s-quadrangular graph. However, $u \not\in O$ by
(iv), $u \not\in R$ by (ii), $u \not\in T$ by (i), and $u \not\in S$
as $e(t,S)=0$. This contradiction implies that $e(t,S)>0$ for all $t
\in T$.

Let $S_1 = \{ s \in S:\ e(s,T) \leq 1 \}$ and let $W=T \cap
N(S-S_1)$. Observe that for every $w \in W$ there exists a vertex $s
\in S-S_1$, such that $w \in N(s)$. Furthermore, there exists a
vertex $w' \in T \cap N(s)-\{w\}$, by the definition of $S_1$. Note
that $w' \in W$, which proves that for every $w \in W$, there is
another vertex, $w' \in W$, such that $N(w) \cap N(w') \not=
\emptyset$. Let $Z=\cup (N(u) \cap N(v):\ u\neq v \in W)$. By (i),
(ii) and (v), $Z \subseteq S-S_1$. By (ii) and (iii), $oc(S,T)\le
e(T,O)$ and, thus, by (vii), $|S|<|T|.$ Observe that $|N(S_1)\cap
T|\le |S_1|$. The last two inequalities and the fact that $e(t,S)>0$
for all $t \in T$ imply $|Z| \leq |S-S_1| < |T-(N(S_1)\cap T)| \leq
|W|$. This is a contradiction to the definition of an
$s$-quadrangular graph. \qed

\2

The next theorem is the main result of this section.

\begin{theorem}
If the degree of every vertex in a connected s-quadrangular graph $G$
is at most 4, then $G$ is hamiltonian.
\end{theorem}
\pf Suppose that $G=(V,E)$ is not hamiltonian. By Theorem
\ref{2-factor}, $G$ has a 2-factor. Let $C_1\cup C_2\cup \cdots \cup
C_m$ be a 2-factor with the minimum number $m\ge 2$ of cycles.
Notice that each cycle $C_i$ is of length at least three.

For a vertex $v$ on $C_i$, $v^+$ is the set of the two neighbors of
$v$ on $C_i$. We will denote the neighbors by $v_1$ and $v_2.$ The
following simple observation is of importance in the rest of the
proof: $$\mbox{If }u,v\mbox{ are vertices of } C_i,C_j,\ i\neq j,
\mbox{ respectively, and }uv\in E,\mbox{ then } e(u^+,v^+)=0.$$
Indeed, if $e(u^+,v^+)>0$, then by deleting $uu^+$ from $C_i$ and
$vv^+$ from $C_j$ and adding an edge between $u^+$ and $v^+$ and the
edge $uv$, we may replace $C_i,C_j$ by just one cycle, which
contradicts minimality of $m.$

We prove that every vertex $u$ which has a neighbor outside its
cycle $C_i$ has degree $4$. Suppose $d_G(u)=3.$ Let $uv \in E$ such
that $v \in C_j$, $j\neq i$. Since $u,v_1$ must have a common
neighbor $z\neq v$, we conclude that $e(u^+,v^+)>0,$ which is
impossible.

Since $G$ is connected, there is a vertex $u \in C_1$ that has a
neighbor outside $C_1$. We know that $d(u)=4$. Apart from the two
vertices in $u^+$, the vertex $u$ is adjacent to two other vertices
$x,y$. Assume that $x,y$ belong to $C_i,C_j$, respectively.
Moreover, without loss of generality, assume that $i\neq 1.$ Since
$x_k$ ($k=1,2$) and $u$ have a common neighbor different from $x$,
$u_1$ and $u_2$, we conclude that $y$ is adjacent to both $x_1$ and
$x_2$. Since $d(y)<5$, we have $|\{u,x_1,x_2,y_1,y_2\}|<5.$ Since
$u,x_1,x_2$ are distinct vertices, without loss of generality, we
may assume that $y_2$ is equal to either $x_1$ or $u.$ If $y_2=u$,
then $y\in u^+$ and $x_1$ is adjacent to a vertex in $u^+$, which is
impossible.

Thus, $y_2=x_1.$ This means that the
vertices $y_1,y,x_1=y_2,x,x_2$ are consecutive vertices of $C_i.$ Recall that $x_2y\in E.$
Suppose that $C_i$ has at least five vertices. Then $y_1,y_2=x_1,x_2$ and $u$ are all distinct neighbors of $y$.
Since $u_1$ and $y$ have a common neighbor different from $u$, the vertex
$u_1$ is adjacent to either $y_1$ or $y_2$ or $x_2$, each possibility implying that
$e(u^+,x^+)+e(u^+,y^+)>0,$ which is impossible.

Thus, $C_i$ has at most four vertices. Suppose $C_i$ has three
vertices: $x,y,x_1.$ Since $u$ and $y$ must have a common neighbor
different from $x$, we have that $y\in x^+$ is adjacent to a vertex
in $u^+$, which is impossible. Thus, $C_i$ has exactly four
vertices: $y$, $x_1=y_2$, $x$ and $x_2=y_1.$ The properties of $u$
and $C_i$ that we have established above can be formulated as the
following general result:

\2

{\bf Claim A:} {\em If a vertex $w$ belonging to a cycle $C_p$
is adjacent to a vertex outside $C_p$, then $w$ is adjacent
to a pair $w'$ and $w''$ of vertices  belonging to a cycle $C_q$
of length four, $q\neq p$, such that $w'$ and $w''$ are not adjacent on $C_q.$}

\2

By Claim A, $x$ and $y$ are adjacent not only to $u$, but also to another vertex $v\notin \{u,u_1,u_2\}$
of $C_1$, and $C_1$ is of length four.

By s-quadrangular property, for the vertex set $S=\{u,x_1,x_2\}$
there must be a set $T$ with at least three vertices such that each
$t\in T$ is a common neighbor of a pair of vertices in $S.$ This
implies that there must be a vertex $z\not\in \{x,y,u_1,u_2,u\}$
adjacent to both $x_1$ and $x_2.$ Thus, $z$ is on $C_k$ with
$k\notin \{1,i\}.$ By Claim A, $x_1$ and $x_2$ are also adjacent to
a vertex $s$ on $C_k$ different from $z_1$ and $z_2$. Continue our
argument with $z_1$ and $z_2$ and similar pairs of vertices, we will
encounter cycles $Z_1,Z_2,Z_3,\ldots$
($Z_j=z_j^1z_j^2z^3_jz^4_jz^1_j$) such that
$\{z^2_jz^1_{j+1},z^4_jz^1_{j+1},z^2_jz^3_{j+1},z^4_jz^3_{j+1}\}\subset
E$ for $j=1,2,3,\ldots.$ Here $Z_1=C_1,Z_2=C_i,Z_3=C_k$.

Since the number of vertices in $G$ is finite, after a while we will encounter a cycle $Z_r$
that we have encountered earlier. We have $Z_r=Z_1$ since for each $m>1$ after we encountered $Z_{m+1}$ we have
established all neighbors of the vertices in $Z_m.$ This implies that
$z^1_1z^2_1z^3_1z^4_1z^1_2z^2_2z^3_2z^4_2\ldots z^1_{r-1}z^2_{r-1}z^3_{r-1}z^4_{r-1}z^1_1$
is a cycle of $G$ consisting of the vertices of the cycles $Z_1,Z_2,\ldots,Z_{r-1}.$ This contradicts minimality
of $m.$
\qed

\2

\2

{\bf Acknowledgments} We are grateful to the referees for remarks
and suggestions that have improved the presentation.

\end{document}